\documentclass[11pt,amssym,twoside]{article}
\usepackage{amssymb}
\newtheorem{thm}{Theorem}[section]
\newtheorem{mthm}[thm]{Main Theorem}

\newtheorem{prop}[thm]{Proposition}
\newtheorem{rem}[thm]{Remark}

\newcommand{\To}{\longrightarrow}

\newcommand{\Z}{\mathbb Z}

\newcommand{\F}{\mathbb F}
\newcommand{\D}{\mathbb D}
\newcommand{\Q}{\mathbb Q}
\newcommand{\LL}{\mathcal {L}}
\newcommand{\N}{\mathcal {N}}
\newcommand{\gal}{\mathcal {G}al}
\newcommand{\out}{\mathcal {O}ut}
\begin{document}

\title{Arithmetic Deformation Theory of Lie Algebras}
\author{Arash Rastegar}


\maketitle
\begin{abstract}
This paper is devoted to deformation theory of graded Lie
algebras over $\Z$ or $\Z_l$ with finite dimensional graded
pieces. Such deformation problems naturally appear in number
theory. In the first part of the paper, we use Schlessinger
criteria for functors on Artin local rings in order to obtain
universal deformation rings for deformations of graded Lie
algebras and their graded representations. In the second part, we
use a version of Schlessinger criteria for functors on the
Artinian category of nilpotent Lie algebras which is formulated
by Pridham, and explore arithmetic deformations using this
technique.
\end{abstract}
\section*{Introduction}

To a hyperbolic smooth curve defined over a number-field one
naturally associates an "anabelian" representation of the
absolute Galois group of the base field landing in outer
automorphism group of the algebraic fundamental group [Gro]. It
would be very fruitful to deform this object, rather than the
abelian version coming from the action of Galois group on the
Tate module of the Jacobian variety. On the other hand, there is
no formalism available for deforming such a representation. This
is why we translate it to the language of graded Lie algebras
which are still non-linear enough to carry more information than
the abelian version. For more details see [Ras1] and [Ras 2].

In this paper, we introduce several deformation problems for
Lie-algebra versions of the above representation. In particular,
we deform the representation of certain Galois-Lie algebra
landing in the graded Lie algebra associated to the weight
filtration on outer automorphism group of the pro-$l$ fundamental
group, and we construct a deformation ring parameterizing all
deformations fixing the mod-$l$ Lie-algebra representation, which
is the typical thing to do in the world of Galois representations.

Organization of the paper is as follows: First, we review methods
of associating Lie algebras to profinite groups. Then, we
introduce some deformation problems for representations landing
in graded Lie algebras. Afterwards, we use the classical
Schlessinger criteria for representability of functors on Artin
local rings for deformation of the above representation. In some
cases, universal deformations exist and in some others, we are
only able to construct a hull which parameterizes all possible
deformations. finally, we use a graded version of Pridham's
adaptation of Schlessinger criteria for functors on finite
dimensional nilpotent Lie algebras.

\section{Lie algebras associated to profinite groups}

Exponential map on the tangent space of an algebraic group
defined over a field $k$ of characteristic zero, gives an
equivalence of categories between nilpotent Lie algebras of
finite dimension over $k$ and unipotent algebraic groups over $k$.
This way, one can associate a Lie algebra to the algebraic
unipotent completion $\Gamma^{alg}(\Q)$ of any profinite group
$\Gamma$.

On the other hand, Mal\v{c}ev defines an equivalence of
categories between nilpotent Lie algebras over $\Q$ and uniquely
dividible nilpotent groups. Inclusion of such groups in nilpotent
groups has a right adjoint $\Gamma\to \Gamma_{\Q}$. For a nipotent
group $\Gamma$, torsion elements form a subgroup $T$ and
$\Gamma_{\Q}= \cup (\Gamma/T)^{1/n}$. In fact we have
$\Gamma_{\Q}=\Gamma^{alg}(\Q)$.

Any nilpotent finite group is a product of its sylow subgroups.
Therefore, the profinite completion $\Gamma^{\wedge}$ factors to
pro-$l$ completions $\Gamma_l^{\wedge}$, each one a compact open
subgroup of the corresponding $\Gamma^{alg}(\Q_l)$ and we have
the following isomorphisms of $l$-adic Lie groups
$$
\mathrm { Lie}(\Gamma_l^{\wedge})=\mathrm {Lie}
(\Gamma^{alg}(\Q_l))=\mathrm {Lie}(\Gamma^{alg}(\Q))\otimes \Q_l.
$$
In fact, the adelic Lie group associated to $\Gamma$ can be
defined as $\mathrm {Lie}(\Gamma^{alg}(\Q))\otimes \mathbb {A}^f$
which is the same as $\prod \mathrm {Lie}(\Gamma_l^{\wedge})$.

Suppose we are given a nilpotent representation of $\Gamma$ on a
finite dimensional vector space $V$ over $k$, which means that
for a filtration $F$ on $V$ respecting the action, the induced
action of $Gr_F(V)$ is trivial. The subgroup
$$
\{\sigma\in GL(V)|\sigma F=F , Gr_F(\sigma)=1\}
$$
is a uniquely divisible group and one obtains a morphism
$$
\mathrm {Lie}(\Gamma_{\Q})\to \{\sigma\in gl(V)| \sigma F=F ,
Gr_F(\sigma)=0\}
$$
which is an equivalence of categories between nilpotent
representations of $\Gamma$ and representations of
$\Gamma^{alg}(\Q)$ and nilpotent representations of the Lie
algebra $\mathrm {Lie}(\Gamma_{\Q})$ over the field $k$. The above
equivalence of categories extends to an equivalence between
linear representations of $\Gamma$ and representations of its
algebraic envelope [Del].

The notion of weighted completion of a group developed by Hain
and Matsumoto generalizes the concept of algebraic unipotent
completion. Suppose that $R$ is an algebraic $k$-group and
$w:\mathbb G_m\to R$ is a central cocharacter. Let $G$ be an
extension of $R$ by a unipotent group $U$ in the category of
algebraic $k$-groups
$$
0\To U\To G\To R\To 0.
$$

The first homology of $U$ is an $R$-module, and therefore an
$\mathbb G_m$-module via $w$, which naturally decomposes to to
direct sum of irreducible representations each isomorphic to a
power of the standard character. We say that our extension is
negatively weighted if only negative powers of the standard
character appear in $H_1(U)$. The weighted completion of $\Gamma$
with respect to the representation $\rho$ with Zariski dense image
$\rho:\Gamma\to R(\Q_l)$ is the universal $\Q_l$-proalgebraic
group $\mathcal G$ which is a negative weighted extension of $R$
by a prounipotent group $\mathcal U$ and a continuous lift of
$\rho$ to $\mathcal G(\Q_l)$ [Hai-Mat]. The Lie algebra of
$\mathcal G(\Q_l)$ is a more sophisticated version of
$\mathrm{Lie}(\Gamma_{\Q})\otimes \Q_l$.

\section{Cohomology theories for graded Lie algebras}

We shall first review cohomology of Lie algebras with the adjoint
representation as coefficients. Let $L$ be a graded $\Z_l$-algebra
and let $C^q(L,L)$ denote the space of all skew-symmetric
$q$-linear forms on a Lie algebra $L$ with values in $L$. Define
the differential
$$
\delta :C^q(L,L)\longrightarrow C^{q+1}(L,L)
$$
where the action of $\delta$ on a skew-symmetric $q$-linear form
$\gamma$ is a skew-symmetric $(q+1)$-linear form which takes
$(l_1,...,l_{q+1})\in L^{q+1}$ to
$$
\sum (-1)^{s+t-1}\gamma ([l_s,l_t],l_1,...,\hat l_s,...,\hat
l_t,...,l_{q+1})+\sum [l_u,\gamma(l_1,...,\hat l_u,...,l_{q+1})]
$$
where the first sum is over $s$ and $t$ with $1\leq s < t \leq
q+1$ and the second sum is over $u$ with $1\leq u\leq q+1$. Then
$\delta^2=0$ and we can define $H^q(L,L)$ to be the cohomology of
the complex $\{C^q(L,L),\delta\}$. If we put $C^m=C^{m+1}(L,L)$
and $H^m=H^{m+1}(L,L)$, then there exists a natural bracket
operation which makes $C=\oplus C^m$ a differential graded
algebra and $H=\oplus H^m$ a graded Lie algebra. Look in
[Fia-Fuc]. If $L$ is $\Z$-graded, $L=\oplus L(m)$, we say $\phi\in
C^q(L,L)(m)$ if for $l_i\in L(g_i)$ we have $\phi(l_1,...,l_q)\in
L(g_1+...+g_q-m)$. Then, there exists a grading induced on the
Lie algebra cohomology $H^q(L,L)=\oplus H^q(L,L)(m)$.

The computational tool used by geometers to study deformations of
Lie algebras is a cohomology theory of $K$-algebras where $K$ is
a field, which is developed by Harrison [Har]. This cohomology
theory is generalized by Barr to algebras over general rings
[Bar]. Here we use a more modern version of the latter introduced
independently by Andre and Quillen which works for general
algebras [Qui].

Let $A$ be a commutative algebra with identity over $\Z_l$ or any
ring $R$ and let $M$ be an $A$-module. By an $n$-long singular
extension of $A$ by $M$ we mean an exact sequence of $A$-modules
$$
0\to M\to M_{n-1}\to ...\to M_1 \to T \to R \to 0
$$
where $T$ is a commutative $\Z_l$-algebra and the final map a
morphism of $\Z_l$-algebras whose kernel has square zero. It is
trivial how to define morphisms and isomorphisms between $n$-long
singular extensions. Baer defines a group structure on these
isomorphism classes [Bar], which defines $H_{Barr}^n(A,M)$ for
$n>1$ and we put $H^1_{Barr}(A,M)=Der(A,M)$. Barr proves that for
a multiplicative subset $S$ of $R$ not containing zero
$$
H_{Barr}^n(A,M)\cong H_{Barr}^n(A_S,M)
$$
for all $n$ and any $A_S$-module $M$. According to this
isomorphism, the cohomology of the algebra $A$ over $\Z_l$ is the
same after tensoring $A$ with $\Q_l$ if $M$ is
$A\otimes\Q_l$-module. Thus one could assume that we are working
with an algebra over a field, and then direct definitions given by
Harrison would serve our computations better. Consider the complex
$$
0\rightarrow Hom(A,M) \rightarrow Hom(S^2A,M) \rightarrow
Hom(A\otimes A\otimes A,M)
$$
where $\psi\in Hom(A,M)$ goes to
$$
d_1\psi:(a,b)\mapsto a\psi(b)-\psi(ab)+b\psi(a)
$$
and $\phi\in Hom(S^2A,M)$ goes to
$$
d_2\phi(a,b,c)\mapsto a\phi(b,c)-\phi(ab,c)+\phi(a,bc)-c\phi(a,b)
$$
The cohomology of this complex defines $H^i_{Harr}(R,M)$ for
$i=1,2$. If $A$ is a local algebra with maximal ideal $m$ and
residue field $k$, the Harrison cohomology
$H^1_{Harr}(A,k)=(m/m^2)'$, which is the space of homomorphisms
$A\to k[t]/t^2$ such that $m$ is the kernel of the composition
$A\to k[t]/t^2\to k$.

Andre-Quillen cohomology is the same as Barr cohomology in low
dimensions and can be described directly in terms of derivations
and extensions. For any morphism of commutative rings $A \to B$
and $B$-module $M$ we denote the $B$-module of $A$-algebra
derivations of $B$ with values in $M$ by $Der_A (B,M)$. Let
$Ext^{inf}_A (B,M)$ denote the $B$-module of infinitesimal
$A$-algebra extensions of $B$ by $M$. The functors $Der$ and
$Ext^{inf}$ have transitivity property. Namely, given morphisms
of commutative rings $A\to B\to C$ and a $C$-module $M$, there is
an exact sequence
$$
0 \longrightarrow Der_B (C,M)\longrightarrow Der_A
(C,M)\longrightarrow Der_A (B,M) \hspace{1 in}$$
$$
\hspace{1 in}\longrightarrow Ext^{inf}_B (C,M)\longrightarrow
Ext^{inf}_A (C,M)\longrightarrow Ext^{inf}_A (B,M).
$$
The two functors $Der$ and $Ext^{inf}$ also satisfy flat
base-change property. Namely, given morphisms $A\to B$ and $A\to
A'$ if $Tor^A_1(A',B)=0$, then there are isomorphisms $Der_{A'}
(A'\otimes_A B,M)\cong Der_A (B,M)$ and $Ext^{inf}_{A'}
(A'\otimes_A B,M)\cong Ext^{inf}_A (B,M)$. Andre-Quillen
cohomology associates $Der_A (B,M)$ and $Ext^{inf}_A (B,M)$ to any
morphism of commutative rings $A \to B$ and $B$-module $M$ as the
first two cohomologies and extends it to higher dimensional
cohomologies such that transitivity and flat base-change extend
in the obvious way.

\section{Several deformation problems}

Let $X$ denote a hyperbolic smooth algebraic curve defined over a
number field $K$. Let $S$ denote the set of bad reduction places
of $X$ together with places above $l$. We shall construct Lie
algebra versions of the pro-$l$ outer representation of the
Galois group
$$
\rho^l_X:Gal(K_S/K) \To Out(\pi_1(\bar X)^{(l)}).
$$
Let $I_l$ denote the decreasing filtration on
$Out(\pi_1(X)^{(l)})$ induced by the central series filtration of
$\pi_1(\bar X)^{(l)}$. By abuse of notation, we also denote the
filtration on $Gal(K_S/K)$ by $I_l$. We get an injection of the
associated graded $\Z_l$-Lie algebras on both sides
$$
\gal (K_S/K) \To \out (\pi_1(\bar X)^{(l)}).
$$
One can also start with the $l$-adic unipotent completion of the
fundamental group and the outer representation of Galois group on
this group.
$$
\rho^{un,l}_X:Gal(K_S/K) \To Out(\pi_1(\bar X)^{un}_{/\Q_l}).
$$
By [Hai-Mat] 8.2 the associated Galois Lie algebra would be the
same as those associated to $I_l$. Let $U_S$ denote the
prounipotent radical of the Zariski closure of the image of
$\rho^{un,l}_X$. The image of $Gal(K_S/K)$ in $Out(\pi_1(\bar
X)^{un}_{/\Q_l})$ is a negatively weighted extension of $\mathbb
G_m$ by $U_S$ with respect to the central cocharacter $w:x\mapsto
x^{-2}$. By [Hai-Mat] 8.4 the weight filtration induces a graded
Lie algebra $\mathcal U_S$ which is isomorphic to $\gal
(K_S/K)\otimes \Q_l$.

There are several deformation problem in this setting which are
interesting. For example, the action of Galois group on unipotent
completion of the fundamental group induces an action of the
Galois group on the corresponding nilpotent $\Q_l$-Lie algebra
$$
\rho^{ni,l}_X:Gal(K_S/K) \To Aut(\mathcal U_S).
$$
which could be deformed. Another possibility is deforming the
following representation
$$
Gal(K_S/K) \To Aut(\pi_1(\bar X)^{un}_{/\Q_l}) \To
Aut(H_1(\mathcal P))
$$
where $\mathcal P$ denote the nilpotent Lie algebra associated to
$\pi_1(\bar X)^{un}_{/\Q_l}$. This time, the Schlessinger
criteria may not help us in finding a universal representation.
There exists also a derivation version, which is a Schlessinger
friendly $\Z_l$-Lie algebras representation
$$
\gal (K_S/K) \To Der(\mathcal P)/Inn(\mathcal P).
$$
One could also deform the following morphism, fixing its mod-$l$
reduction
$$
\gal (K_S/K) \To \out (\pi_1(\bar X)^{(l)}).
$$

\section{Deformations of local graded Lie algebras}

In this section, we are only concerned with deformations of Lie
algebras and leave deformation of their representations for the
next section. We are interested in deforming the coefficient ring
of graded Lie algebras over $\Z_l$ of the form $L=Gr^{\bullet}_I
\widetilde {Out}(\pi^l_1(X))$ and then deforming representations
of Galois graded Lie algebra
$$
Gr^{\bullet}_{X,l} Gal(\bar {K}/K)\to Gr^{\bullet}_I \widetilde
{Out}(\pi^l_1(X)).
$$
One can reduce the coefficient ring $\Z_l$ modulo $l$ and get a
graded Lie algebra $\bar L$ over $\F_l$ and a representation
$$
Gr^{\bullet}_{X,l} Gal(\bar {K}/K)\to \bar {L} .
$$
We look for liftings of this representation which is landing in
$\bar L$ among representations landing in graded Lie algebras
over Artin local rings $A$ of the form $L=\oplus_i L^i$ where
$L^i$ is a finitely generated $A$-module for positive $i$.

Let $k$ be a finite field of characteristic $p$ and let $\Lambda$
be any complete Noetherian local ring. For example $\Lambda$ can
be $W(k)$, the ring of Witt vectors of $k$, or $O$, the ring of
integers of any local field with quotient field $K_{\wp}$ and
residue field $k$. Let $C$ denote the category of Artinian local
$\Lambda$-algebras with residue field $k$. A covariant functor
from $C$ to $Sets$ is called pro-representable if it has the form
$$
F(A)\cong Hom_{\Lambda}(R_{univ},A)\hspace {0.5 in} A\in C
$$
where $R_{univ}$ is a complete local $\Lambda$-algebra with
maximal ideal $m_{univ}$ such that $R_{univ}/m_{univ}^n$ is in
$C$ for all $n$.

There are a number of deformation functors related to our problem.
For a local ring $A\in C$ with maximal ideal $m$, the set of
deformations of $\bar L$ to $A$ is denoted by $D_c(\bar L,A)$ and
is defined to be the set of isomorphism classes of graded Lie
algebras $L/A$ of the above form which reduce to $\bar L$ modulo
$m$. In this notation $c$ stands for coefficients, since we are
only deforming coefficients not the Lie algebra structure. The
functor $D(\bar{L})$ as defined above is not a pro-representable
functor. As we will see, there exists a "hull" for this functor
(Schlessinger's terminology [Sch]) parameterizing all possible
deformations. The idea of deforming the Lie structure of Lie
algebras has been extensively used by geometers. For example,
Fialowski studied this problem in double characteristic zero case
[Fia]. In this paper, we are interested in double characteristic
$(0,l)$-version.

Now we make an assumption for further constructions. Assume
$H^2(L,L)(m)$ is finite dimensional for all $m$ and consider the
algebra $\D_1=\Z_l\oplus\bigoplus_m H^2(L,L)(m)'$ where $'$ means
the dual over $\Z_l$. Fix a graded homomorphism of degree zero
$$
\mu :H^2(L,L)\longrightarrow C^2(L,L)
$$
which takes any cohomology class to a cocycle representing this
class. Now define a Lie algebra structure on
$$
\D_1\otimes L=L\oplus Hom^0(H^2(L,L),L)
$$
where $Hom^0$ means degree zero graded homomorphisms, by the
following bracket
$$
[(l_1,\phi_1),(l_2,\phi_2)]:=([l_1,l_2],\psi )
$$
where $\psi(\alpha )=\mu(\alpha)(l_1,l_2)+[\phi_1(\alpha),l_2]
+[\phi_2(\alpha),l_1]$. The Jacobi identity is implied by
$\delta\mu(\alpha)=0$. It is clear that this in an infinitesimal
deformation of $L$ and it can be shown that, up to an isomorphism,
this deformation does not depend on the choice of $\mu$. We shall
denote this deformation by $\eta_L$ after Fialowski and Fuchs
[Fia-Fuc].
\begin{prop} Any infinitesimal deformation of $\bar {L}$ to a
finite dimensional local ring $A$ is induced by pushing forward
$\eta_L$ by a unique morphism
$$
\phi:\Z_l\oplus\bigoplus_m H^2(L,L)(m)'\longrightarrow A.
$$
\end{prop}
\textbf {Proof.} This is the double characteristic version of
proposition 1.8 in [Fia-Fuc]. $\square$

Note that, in our case $H^2(L,L)$ is not finite dimensional. This
is why we restrict our deformations to the space of graded
deformations. Since $H^2(L,L)(0)$ is the tangent space of the
space of graded deformations, and the grade zero piece
$H^2(L,L)(0)$ is finite dimensional, the following version is more
appropriate:
\begin{prop} Any infinitesimal graded deformation of $\bar {L}$ to a
finite dimensional local ring $A$ is induced by pushing forward
$\eta^0_L$ by a unique morphism
$$
\phi:\Z_l\oplus H^2(L,L)(0)'\longrightarrow A.
$$
where $\eta^0_L$ denotes the restriction of $\eta_L$ to
$\Z_l\oplus H^2(L,L)(0)'$.
\end{prop}

Let $A$ be a small extension of $\mathbb {F}_l$ and $L$ be a
graded deformation of $\bar{L}$ over the base $A$. The
deformation space $D(\bar {L},A)$ can be identified with
$H^2(L,L)(0)$ which is finite dimensional. Therefore, by
Schlessinger criteria, in the subcategory $C'$ of $C$ consisting
of local algebras with $m^2=0$ for the maximal ideal $m$, the
functor $D(\bar {L},A)$ is pro-representable. This means that
there exists a unique map
$$
\Z_l\oplus H^2(L,L)(0)'\to R_{univ}
$$
inducing the universal infinitesimal graded deformation.

\subsection{Obstructions to deformations}

Let $A$ be an object in the category $C$ and $L\in Def(\bar
{L},A)$. The pair $(A,L)$ defines a morphism of functors $\theta
:Mor(A,B)\to Def(\bar {L},B)$. We say that $(A,L)$ is universal if
$\theta$ is an isomorphism for any choice of $B$. We say that
$(A,L)$ is miniversal if $\theta$ is always surjective, and gives
an isomorphism for $B=k[\varepsilon]/\varepsilon^2$. We intend to
construct a miniversal deformation of $\bar {L}$.

Consider a graded deformation with base in a local algebra $A$
with residue field $k=\F_l$. One can define a map
$$
\Phi_A:Ext^{inf}_{\Z_l}(A,k)\longrightarrow H^3(\bar {L},\bar
{L}).
$$
Indeed, choose an extension $0\to k\to B\to A \to 0$ corresponding
to an element in $Ext^{inf}_{\Z_l}(A,k)$. Consider the $B$-linear
skew-symmetric operation $\{ .,.\}$ on $\bar {L}\otimes_k B$
commuting with $[.,.]$ on $\bar {L}\otimes A$ defined by $\{
l,l_1\}=[l,\bar {l}_1]$ for $l$ in the kernel of $\bar
{L}\otimes_k B\to \bar {L}\otimes_k A$ which can be identified by
$\bar {L}$. Here $\bar {l}_l$ is the image of $l_1$ under the
projection map $B\to k$ tensored with $\bar L$ whose kernel is the
inverse image of the maximal ideal of $A$. The Jacobi expression
induces a multilinear skew-symmetric form on $\bar L$ which could
be regarded as a closed element in $C^3(\bar {L},\bar {L})$. The
image in $H^3(\bar {L},\bar {L})$ is independent of the choices
made.
\begin{thm}(Fialowski)
One can deform the Lie algebra structure on $L\otimes_k A$ to
$L\otimes_k B$ if and only if the image of the above extension
vanishes under the morphism $Ext^{inf}_{\Z_l}(A,k)\to H^3(\bar
{L},\bar {L})$.
\end{thm}
\textbf {Proof.} The proof presented in [Fia] and [Fia-Fuc] works
for algebras over fields of finite characteristic. $\square$

Using the above criteria for extending deformations, one can
follow the methods of Fialowski and Fuchs to introduce a
miniversal deformation for $\bar L$.
\begin{prop}
Given a local commutative algebra $A$ over $\Z_l$ there exists a
universal extension
$$
0 \longrightarrow Ext^{inf}_{\Z_l}(A,k)'\longrightarrow C
\longrightarrow A \longrightarrow 0
$$
among all extensions of $A$ with modules $M$ over $A$ with $mM=0$
where $m$ is the maximal ideal of $A$.
\end{prop}
This is proposition 2.6 in [Fia-Fuc]. Consider the canonical split
extension
$$
0\longrightarrow H^2(L,L)(0)'\longrightarrow \D_1\longrightarrow
k\longrightarrow 0.
$$
We will initiate an inductive construction of $\D_k$ such that
$$
0\longrightarrow Ext^{inf}_{\Z_l}(\D_k,k)'\longrightarrow \bar
{\D}_{k+1}\longrightarrow \D_k\longrightarrow 0.
$$
together with a deformation $\eta_k$ of $L$ to the base $\D_k$.
For $\eta_1$ take $\eta_L$, and assume $\D_i$ and $\eta_k$ is
constructed for $i\leq k$. Given a local commutative algebra
$\D_k$ with maximal ideal $m$, there exists a unique universal
extension for all extensions of $\D_k$ by $\D_k$-modules $M$ with
$mM=0$ of the following form
$$
0\longrightarrow Ext^{inf}_{\Z_l}(\D_k,k)'\longrightarrow C
\longrightarrow \D_k\longrightarrow 0.
$$
associated to the cocycle $f_k:S^2(\D_k)\to
Ext^{inf}_{\Z_l}(\D_k,k)'$ which is dual to the homomorphism
$$
\mu: Ext^{inf}_{\Z_l}(\D_k,k)\longrightarrow S^2(\D_k)'
$$
which takes a cohomology class to a cocycle from the same class.
The obstruction to extend $\eta_k$ lives in
$Ext^{inf}_{\Z_l}(\D_k,k)'\otimes H^3(\bar {L},\bar {L})$.
Consider the composition of the associated dual map
$$
\Phi'_k:H^3(\bar {L},\bar {L})'\longrightarrow
Ext^{inf}_{\Z_l}(\D_k,k)'.
$$
with $Ext^{inf}_{\Z_l}(\D_k,k)'\to  C$ and define $\D_{k+1}$ to
be the cokernel of this map. We get the following exact sequence
$$
0\longrightarrow (ker\Phi_k)'\longrightarrow \D_{k+1}
\longrightarrow \D_k\longrightarrow 0.
$$
We can extend $\eta_k$ to $\eta_{k+1}$. Now, taking a projective
limit of $\D_k$ we get a base and a formal deformation of $\bar
{L}$.
\begin{thm}
Let $\D$ denote the projective limit $\underline {\lim} \D_k$
which is a $\Z_l$-module. One can deform $\bar L$ uniquely to a
graded Lie algebra with base $\D$ which is miniveral among all
deformations of $\bar {L}$ to local algebras over $\Z_l$.
\end{thm}
\textbf {Proof.} This is the double characteristic version of
theorem 4.5 in [Fia-Fuc]. The same proof works here because
theorems 11 and 18 in [Har] which are used in the arguments of
Fialowski and Fuchs work for algebras over any perfect field.
$\square$

\begin{prop} (Fialowski-Fuchs) The base of the minversal
deformation of $\bar {L}$ is the zero locus of a formal map
$H^2(\bar {L},\bar {L})(0)\to H^3(\bar {L},\bar {L})(0)$.
\end{prop}
\textbf {Proof.} This is the graded version of proposition 7.2 in
[Fia-Fuc].$\square$

\subsection{Deformations of graded Lie algebra representations}

In the previous section we discussed deformation theory of the
mod $l$ reduction of the Lie algebra $Gr^{\bullet}_I \widetilde
{Out}(\pi^l_1(X))$. We shall mention the following
\begin{thm}  The cohomology
groups $H^i(Gr^{\bullet}_I \widetilde
{Out}(\pi^l_1(X)),Gr^{\bullet}_I \widetilde
{Out}(\pi^l_1(X)))(0)$ are finite dimensional for all non-negative
integer $i$.
\end{thm}
\textbf{Proof.} By a theorem of Labute $Gr^{\bullet}_I \pi^l_1(X)$
is quotient of a finitely generated free Lie algebra with
finitely generated module of relations [Lab]. Therefore, the
cohomology groups $H^i(Gr^{\bullet}_I \pi^l_1(X),Gr^{\bullet}_I
\pi^l_1(X))(0)$ are finite dimensional. Finite dimensionality of
the cohomology of $Gr^{\bullet}_I \widetilde {Out}(\pi^l_1(X))$
follows from proposition 1.3. $\square$

We are interested in deforming the following graded
representation of the Galois graded Lie algebra
$$
\rho:Gr^{\bullet}_{X,l} Gal(\bar {K}/K)\to Gr^{\bullet}_I
\widetilde {Out}(\pi^l_1(X))
$$
among all graded representations which modulo $l$ reduce to the
graded representation
$$
\bar {\rho}:Gr^{\bullet}_{X,l} Gal(\bar {K}/K)\to \bar {L}
$$
where the Lie algebra $\bar L$ over $\F_l$ is the mod-$l$
reduction of $Gr^{\bullet}_I \widetilde {Out}(\pi^l_1(X))$. There
are suggestions from the classical deformation theory of Galois
representations on how to get a representable deformation functor.
Let $D(\bar {\rho},A)$ denote the set of isomorphism classes of
Galois graded Lie algebra representations to graded Lie algebras
$L/A$ of the above form which reduce to $\bar {\rho}$ modulo $m$.
The first ingredient we need in order to prove representability of
$D(\bar {\rho})$ is finite dimensionality of the tangent space of
the functor. The tangent space of the deformation functor $D(\bar
{\rho})$ for an object $A\in C$ is canonically isomorphic to
$$
H^1(Gr^{\bullet}_{X,l} Gal(\bar {K}/K),Ad\circ \overline {\rho})
$$
where the Lie algebra module is given by the composition of $\bar
{\rho}$ with the adjoint representation of $Gr^{\bullet}_I
\widetilde {Out}(\pi^l_1(X/K))$. To get finite dimensionality, we
restrict ourselves to the graded deformations of the graded
representation $\bar {\rho}$.
\begin{thm}
$H^1(Gr^{\bullet}_{X,l} Gal(\bar {K}/K),Ad\circ \overline
{\rho})(0)$ is finite dimensional.
\end{thm}
\textbf {Proof.} The Galois-Lie representation $Gr^{\bullet}_I
Gal(\bar {K}/K)\to Gr^{\bullet}_I \widetilde {Out}(\pi^l_1(X))$
is an injection. Derivation inducing cohomology commutes with
inclusion of Lie algebras. Therefore $H^1(Gr^{\bullet}_{X,l}
Gal(\bar {K}/K),Ad\circ \overline {\rho})(0)$ injects in
$H^1(Gr^{\bullet}_I \widetilde {Out}(\pi^l_1(X)),Ad\circ
\overline {\rho})(0)$ which is finite dimensional by previous
theorem. $\square$

For a surjective mapping $A_1\to A_0$ of Artinian local rings in
$C$ such that the kernel $I\subset A_1$ satisfies $I.m_1=0$ and
given any deformation $\rho_0$ of $\bar {\rho}$ to $L\otimes A_0$
one can associate a canonical obstruction class in
$$
H^2(Gr^{\bullet}_{X,l} Gal(\bar {K}/K),I\otimes Ad\circ \bar
{\rho})\cong H^2(Gr^{\bullet}_{X,l} Gal(\bar {K}/K),Ad\circ \bar
{\rho})\otimes I
$$
which vanishes if and only if $\rho_0$ can be extended to a
deformation with coefficients $A_1$. Therefore, vanishing results
on second cohomology are important.
\begin{thm}
Suppose $Gr^{\bullet}_{X,l} Gal(\bar {K}/K)$ is a free Lie algebra
over $\Z_l$, then the Galois cohomology $H^2(Gr^{\bullet}_{X,l}
Gal(\bar {K}/K),Ad\circ \bar {\rho})$ vanishes.
\end{thm}
\textbf {Proof.} The free Lie algebra $G=Gr^{\bullet}_{X,l}
Gal(\bar {K}/K)$ is rigid, and therefore has trivial infinitesimal
deformations. Thus, we get vanishing of its second cohomology:
$$
H^2(Gr^{\bullet}_{X,l} Gal(\bar {K}/K),Gr^{\bullet}_{X,l} Gal(\bar
{K}/K))=0.
$$
The injection of $G$ inside $L=Gr^{\bullet}_I \widetilde
{Out}(\pi^l_1(X))$ as Lie-algebras over $\Z_l$ implies that, the
cohomology group $H^2(Gr^{\bullet}_{X,l} Gal(\bar {K}/K),Ad\circ
\rho)$ vanishes again by freeness of $G$. Let $\bar {G}$ denote
the reduction modulo $l$ of $G$ which is a free Lie algebra over
$\F_l$. The cohomology $H^2(G,\bar {G})$ is the mod-$l$ reduction
of $H^2(G,G)$, hence it also vanishes. So does the cohomology
$H^2(Gr^{\bullet}_{X,l} Gal(\bar {K}/K),Ad\circ \bar {\rho})$ by
similar reasoning.$\square$

We have obtained conceptual conditions implying $H_4$ of [Sch].
What we have proved can be summarized as follows.

\begin{mthm} Suppose that $Gr^{\bullet}_{X,l} Gal(\bar {K}/K)$ is a
free Lie algebra over $\Z_l$. There exists a universal deformation
ring $R_{univ}=R(X,K,l)$ and a universal deformation of the
representation $\bar {\rho}$
$$
\rho^{univ}:Gr^{\bullet}_{X,l} Gal(\bar {K}/K)\longrightarrow
Gr^{\bullet}_I \widetilde {Out}(\pi^l_1(X))\otimes R_{univ}
$$
which is unique in the usual sense. If $Gr^{\bullet}_{X,l}
Gal(\bar {K}/K)$ is not free, then a mini-versal deformation
exists which is universal among infinitesimal deformations of
$\bar {\rho}$.
\end{mthm}
\begin{rem}
Note that, freeness of $Gr^{\bullet}_{X,l} Gal(\bar {K}/K)$ in the
special case of $K=\Q$ where filtration comes from punctured
projective curve $X=\mathbb {P}^1-\{ 0,1,\infty\}$ or a punctured
elliptic curve $X=E-\{ 0\}$ is implied by Deligne's conjecture.
\end{rem}

As in the classical case, the Lie algebra structure on
$Ad\circ\bar{\rho}$ induces a graded Lie algebra structure on the
cohomology $H^*(Gr^{\bullet}_{X,l} Gal(\bar {K}/K),Ad\circ \bar
{\rho})$ via cup-product, and in particular, a symmetric bilinear
pairing
$$
H^1(Gr^{\bullet}_{X,l} Gal(\bar {K}/K),Ad\circ \bar {\rho})\times
H^1(Gr^{\bullet}_{X,l} Gal(\bar {K}/K),Ad\circ \bar {\rho})
\hspace {1 in}
$$
$$
\hspace {3 in}\longrightarrow H^2(Gr^{\bullet}_{X,l} Gal(\bar
{K}/K),Ad\circ \bar {\rho})
$$
which gives the quadratic relations satisfied by the minimal set
of formal parameters of $R_{univ}/lR_{univ}$ for characteristic
$l$ different from $2$.

\section{Functors on nilpotent graded Lie algebras}

In this section, we review Pridham's nilpotent Lie algebra version
of Schlessinger criteria [Pri]. The only change we impose is to
consider finitely generated graded nilpotent Lie algebras with
finite dimensional graded pieces, instead of finite dimensional
nilpotent Lie algebras.

Fix a field $k$ and let $\N_k$ denote the category of finitely
generated NGLAs (nilpotent graded Lie algebras) with finite
dimensional graded pieces, and $\widehat {\N_k}$ denote the
category of pro-NGLAs with finite dimensional graded pieces which
are finite dimensional in the sense that $\dim L/[L,L]<\infty$.
Given $\LL \in \widehat {\N_k}$ define $\N_{\LL,k}$ to be the
category of pairs $\{ N\in \N_k, \phi :\LL\to N\}$ and $\widehat
{\N_{\LL,k}}$ to be the category of pairs $\{ N\in \widehat
{\N_k}, \phi :\LL\to N\}$.

All functors on $\N_{\LL,k}$ should take the $0$ object to a one
point set. for a functor $F:\N_{\LL,k}\to \textrm {Set}$, define
$\hat F:\widehat {\N_{\LL,k}}\to \textrm {Set}$ by
$$
\hat F(L)=\lim_{\leftarrow} F(L/\Gamma_n(L)),
$$
where $\Gamma_n(L)$ is the $n$-th term in the central series of
$L$. Then for $h_L:\N_{\LL,k}\to \textrm {Set}$ defined by $N\to
Hom(L,N)$ we have an isomorphism
$$
\hat F(L) \longrightarrow Hom(h_L,F)
$$
which can be used to define the notion of a pro-representable
functor.

A morphism $p\in N\to M$ in $\N_{\LL,k}$ is called a small
section if it is surjective with a principal ideal kernel $(t)$
such that $[N,(t)]=(0)$.

Given $F:\N_{\LL,k}\to \textrm {Set}$, and morphisms $N'\to N$
and $N''\to N$ in $\N_{\LL,k}$, consider the map
$$
F(N'\times_N N'')\To F(N')\times_{F(N)} F(N'').
$$
Then, by the Lie algebra analogue of the Schlessinger theorem $F$
has a hull if and only if it satisfies the following properties

(H1) The above map is surjective whenever $N''\to N$ is a small
section.

(H2) The above map is bijective when $N=0$ and $N''=L(\epsilon )$.

(H3) $\dim_k (t_F< \infty$. \\ $F$ is pro-representable if and
only if it satisfies the following additional property

(H4) The above map is an isomorphism for any small extension
$N''\to N$. \\ Note that, in case we are considering graded
deformations of graded Lie algebras, only the zero grade piece of
the cohomology representing the tangent space shall be checked to
be finite dimensional.


\section{Deformation of nilpotent graded Lie algebras}

Let us concentrate on deforming the graded Lie algebra
representation
$$
\gal (K_S/K) \To \out (\pi_1(\bar X)^{(l)}).
$$
We could fix the mod-$l$ representation, or fix restriction of
this representation to decomposition Lie algebra $\mathcal D_p$
at prime $p$, which is induced by the same filtration as
$Gal(K_S/K)$ on the decomposition group. For each prime $p$ of
$K$ we get a map
$$\mathcal D_p \To \gal (K_S/K)\To \out (\pi_1(\bar X)^{(l)}).
$$
\begin{thm}
For a graded $\Z_l$-Lie algebra $L$, let $D(L)$ be the set of
representations of $\gal(K_S/K)$ to $L$ which reduce to
$$
\bar {\rho}: \gal (K_S/K) \To \out (\pi_1(\bar X)^{(l)})/l\out
(\pi_1(\bar X)^{(l)})
$$
after reduction modulo $l$. Assume that $\gal(K_S/K)$ is a free
$\Z_l$-Lie algebra. Then, there exists a universal deformation
graded $\Z_l$-Lie algebra $L_{univ}$ and a universal
representation
$$
\gal (K_S/K) \To L_{univ}
$$
representing the functor $D$. In case $\gal(K_S/K)$ is not free,
then one can find a hull for the functor $D$.
\end{thm}
\textbf {Proof.} For free $\gal (K_S/K)$ by theorem 4.9, we have
$$
H^2(\gal (K_S/K), Ad\circ \bar{\rho} )=0
$$
which implies that $D$ is pro-representable. In case $\gal
(K_S/K)$ is not free, we have constructed a miniversal
deformation Lie algebra for another functor , which implies that
the first three Schlessinger criteria hold. By a similar argument
one could prove that there exists a hull for $D$. Note that $\out
(\pi_1(\bar X)^{(l)})$ is pronilpotent, and for deformation of
such an object one should deform the truncated object and then
take a limit to obtain a universal object in pro-NGLAs. $\Box$

Note that, by a conjecture of Deligne, there exists a graded Lie
algebra over $\Z$ which gives rise to all $\gal (K_S/K)$ for
different primes $l$, after tensoring with $\Z_l$. For deformation
of the representations of this $\Z$-Lie algebra, one can not use
representability for functors on Artin local ring, and the
nilpotent Lie algebra deformations become crucial.

\section*{Acknowledgements}
I would like to thank M. Kontsevich, O. Gabber, Y. Soibelman for
enjoyable conversations. Also, I wish to thank abdus salam
international center for theoretical physics and institute des
hautes \`{e}tudes scientifiques for their warm hospitality during
which this work was prepared.

Sharif University of Technology, e-mail: rastegar@sharif.edu

\end{document}